\documentclass[12pt]{article}
\usepackage{amssymb}
\usepackage{amsmath}

\newcommand{\be}{\begin{eqnarray}}
\newcommand{\ee}{\end{eqnarray}}
\newcommand{\gf}[2]{f_{#2}^{[#1]}}
\newcommand{\igf}[2]{g_{#2}^{[#1]}}
\newcommand{\az}{A\hspace{-0.28cm}Z}

\newtheorem{theo}{Theorem}

\newtheorem{defn}{Definition}

\title{The Number of [Old-Time] Basketball games with Final Score n:n where the Home Team was never losing but also never ahead by more than w Points}
\author{Arvind Ayyer\\
\small Department of Physics\\
\small 136 Frelinghuysen Rd\\
\small Piscataway, NJ 08854.\\
\small \texttt{ayyer@physics.rutgers.edu}\\
\\
Doron Zeilberger\\
\small Department of Mathematics\\
\small 110 Frelinghuysen Rd\\
\small Piscataway, NJ 08854.\\
\small \texttt{zeilberg@math.rutgers.edu}\\
}

\begin{document}

\maketitle

\begin{abstract}
We show that the generating function (in $n$) for the number of walks on the square lattice with steps $(1,1), (1,-1), (2,2)$ and $(2,-2)$ from $(0,0)$ to $(2n,0)$ in the region $0 \leq y \leq w$ satisfies a very special fifth order \emph{nonlinear} recurrence relation in $w$ that implies both its numerator and denominator satisfy a \emph{linear} recurrence relation.
\end{abstract}

\section{Introduction}
We consider walks in the two-dimensional square lattice with steps (1,1), (1,-1), (2,2) and (2,-2). We assign a weight $\sqrt{z}$ for a unit distance along the $x$-axis. We constrain them to lie in the region defined by $y \geq 0$ and $y \leq w$. The motivation for considering such walks is the modelling of polymers forced to lie between plates separated by a small distance. 

One would then like to calculate various combinatorial quantities. In principle, one hopes to count all possible configurations of the polymer modelling it as a self-avoiding walk [WSCM, MTW]. Since this is a tough nut to crack, one simplifying approach is to treat the polymer as a \emph{directed walk}.

Studies of this kind have been done in the literature with simpler steps such as Dyck paths ($(1,1)$ and $(1,-1)$) which we review in the next section. See, for example, [DR, BORW]. For further developments on the subject, see [R] and the references therein.

Even though the motivation came from Physics [BORW], it later occured to us that this is the number of basketball games (post-1896 and pre-1961, when the three-pointer did not exist) in which the home team always leads the visitor by at most $w$ points ending in a tie!

\section{$\mathrm{\acute{E}tude}$ - Soccer Games}
As a warm-up to the study of basketball games, let us consider soccer games with the same condition [BORW]. These are exactly Dyck walks on the square lattice restricted to $0 \leq y \leq w$ starting at the origin and ending on the $x$-axis. As is usual, we assign a weight $\sqrt{z}$ for both steps.

Let $C_w(z)$ be the generating function for such a walk. And $D_w(z)$ be the generating function for an irreducible walk. That is, one which does not touch the $x$-axis in the interior of the walk. A general walk is either the null walk or is composed of an irreducible walk followed by a smaller such walk. Thus,
\be
C_w = 1 + D_w C_w.
\ee

And an irreducible walk starts with the $(1,1)$ step and ends with the $(1,-1)$ step with an arbitrary walk in between whose width is $w-1$.
\be
D_w & = & \sqrt{z} C_{w-1} \sqrt{z} \cr
& = & z C_{w-1},
\ee
\noindent
which implies
\be
C_w = \frac{1}{1- z C_{w-1}}. \label{Cw}
\ee

This leads to a nice continued fraction expression for $C_w$, which has the distinct aroma of Tchebyshev! Notice that $C_0 = 1$ and thus, $C_1 = 1/(1-z)$. Then
\be
C_w = \frac{1}{1-} \underbrace{\frac{z}{1-} \cdots \frac{z}{1-}}_{w-2 \, \mathrm{terms}} \frac{z}{1-z} \; \mathrm{for} \; w \geq 2 
\ee
\noindent

(\ref{Cw}) is a patently nonlinear recurrence for the generating function. But it does lead to a \emph{linear} recurrence for the numerator and denominator of $C_w$. This can be seen by setting $C_w = \frac{P_w}{Q_w}$. It is easily seen (do it!) that the linear recurrence relations
\be
P_w & = & Q_{w-1} \\
Q_w & = & Q_{w-1} - z Q_{w-2}
\ee
\noindent
with suitable initial conditions gives rise to $C_w$. Notice that these are recurrences with \emph{constant coefficients} in $w$ but, of course, not in $z$. This explains the relationship of the denominators with Tchebyshev polynomials of the first kind - $T_n(z)$ which satisfies a very similar second order recurrence relation in $n$ with constant coefficients, viz.
\be
T_n(z) = 2zT_{n-1}(z) - T_{n-2}(z).
\ee

As an aside, note that if  $w=2$, the number of walks ending at $(n,0)$ give rise to the Fibonacci numbers and if $w=\infty$, the  Catalan numbers [St].

\section{{\it The Main Result} - Basketball Games}

\begin{defn} An \emph{$[ij]$ walk} is a walk that starts at the line $y=i$ and ends at the line $y=j$.
\end{defn}

\begin{defn} An \emph{irreducible} $[ij]$ walk is an $[ij]$ walk that touches the minimum of $i$ and $j$ only at the corresponding endpoint.
\end{defn}

We will need various kinds of generating functions in the proof. Let $\gf{ij}{w}(z)$ denote the generating function of the $[ij]$ walk with width $w$. And let $\igf{ij}{w}(z)$ denote the generating function for the corresponding \emph{irreducible} version of the walk. Note that, at the end of the day, we need a recurrence relation for $F_w := \gf{00}{w}$.

\begin{theo} Let $F_w$ be defined as above. Then it satisfies the following recurrence relation.
\be
& F_w = 1- z F_w + 2 z F_w F_{w-1} + 2 z^2 F_w F_{w-1} F_{w-2}\cr
& -(z^3+z^4)F_w F_{w-1} F_{w-2} F_{w-3} + z^5 F_w F_{w-1} F_{w-2} F_{w-3} F_{w-4} \label{Fw}
\ee
\end{theo}

To prove this, we first write down a set of equations relating different generating functions and then try to solve for $F_w$. First off, a $[00]$ walk is either the empty walk or it is composed of an irreducible $[00]$ walk  followed by a smaller $[00]$ walk.

\be
\gf{00}{w} & = & 1 + \gf{00}{w} \igf{00}{w} \label{f00}
\ee

Next, a $[01]$ walk is always uniquely composed of an arbitrary $[00]$ walk followed by an irreducible $[01]$ walk. Similarly, a $[10]$ walk is uniquely composed of an irreducible $[10]$ walk followed by an arbitrary $[00]$ walk. 

\be
\gf{10}{w} & = & \igf{10}{w} \gf{00}{w} \label{f01}\\
\gf{01}{w} & = & \igf{01}{w} \gf{00}{w} \label{f10}
\ee

A $[11]$ walk either never goes below the first level, in which case it is simply the same as a $[00]$ walk with width $w-1$, or if it does, it is composed of an irreducible $[10]$ walk followed by an arbitrary $[01]$ walk. 

\be
\gf{11}{w} & = & \gf{00}{w-1} + \igf{10}{w} \gf{01}{w} \label{f11}
\ee

Now, we go on to describe the irreducible walks. Since we have a finite width, we will describe them in terms of generating functions for lower widths. In each case, we have to consider different cases for the starting step and the ending step. First, an irreducible $[00]$ walk can begin with either the $(1,1)$ or $(2,2)$ step and end with either the $(1,-1)$ or $(2,-2)$ step. If the walk starts with $(1,1)$ and ends with $(1,-1)$, then there could be an arbitrary $[00]$ walk with width $w-1$ in between. If the walk starts with $(1,1)$ and ends with $(2,-2)$, there has to be an arbitrary $[01]$ walk with width $w-1$ in between. If the walk starts with $(2,2)$ and ends with $(1,-1)$, there has to be an arbitrary $[10]$ walk with width $w-1$ in between. And finally, if the walk starts with $(2,2)$ and ends with $(2,-2)$, there is a $[11]$ walk with width $w-1$ in between.

\be
\igf{00}{w} & = & z \gf{00}{w-1} + z^{3/2} \gf{01}{w-1} + z^{3/2} \gf{10}{w-1} + z^2 \gf{11}{w-1} \label{g00}
\ee

For an irreducible $[01]$ walk, we just need to consider the starting steps. If it starts with $(1,1)$, the remainder is an arbitrary $[00]$ walk with width $w-1$. If it starts with $(2,2)$, the remainder is again an arbitrary $[10]$ walk with width $w-1$. A very similar argument on the ending step yields the equation for an irreducible $[10]$ walk.

\be
\igf{01}{w} & = & z^{1/2} \gf{00}{w-1} + z \gf{10}{w-1} \label{g01}\\
\igf{10}{w} & = & z^{1/2} \gf{00}{w-1} + z \gf{01}{w-1} \label{g10}
\ee

First we eliminate the irreducible generating functions using equations (\ref{g00}), (\ref{g01}) and (\ref{g10}). Then equations (\ref{f00}), (\ref{f01}),  (\ref{f10}) and  (\ref{f11}) become 
\be
\gf{00}{w} & = & 1 + \gf{00}{w} (z \gf{00}{w-1} + z^{3/2} \gf{01}{w-1} + z^{3/2} \gf{10}{w-1} + z^2 \gf{11}{w-1}) \label{2f00}\\
\gf{01}{w} & = & \gf{00}{w}(z^{1/2} \gf{00}{w-1} + z \gf{01}{w-1}) \label{2f01} \\
\gf{10}{w} & = & \gf{00}{w}(z^{1/2} \gf{00}{w-1} + z \gf{10}{w-1}) \label{2f10} \\
\gf{11}{w} & = & \gf{00}{w-1} + \gf{01}{w} (z^{1/2} \gf{00}{w-1} + z \gf{01}{w-1}) \label{2f11}.
\ee

We clean up our notation now. Let $F_w := \gf{00}{w},G_w := \gf{01}{w},H_w := \gf{10}{w},J_w := \gf{11}{w}$. Then
\be
F_w & = & 1 + F_w (z F_{w-1} + z^{3/2} G_{w-1} + z^{3/2} H_{w-1} + z^2 J_{w-1}) \label{F} \\
G_w & = & F_w (z^{1/2} F_{w-1} + z H_{w-1}) \label{G} \\
H_w & = & F_w (z^{1/2} F_{w-1} + z G_{w-1}) \label{H} \\
J_w & = & F_{w-1} + G_w (z^{1/2} F_{w-1} + z G_{w-1}) \label{J}.
\ee

Using (\ref{G}) and (\ref{H}),
\be
G_w - H_w = z F_w (H_{w-1} - G_{w-1}). \label{GH}
\ee

But notice that $G_0 = H_0 = 0$ by definition. Therefore, inductively, $G_w = H_w$. Thus,
\be
G_w & = & F_w (z^{1/2} F_{w-1} + z G_{w-1}) \label{G2}.
\ee

Now we eliminate everything in (\ref{F}) in the form of $G_w$ and $F_w$ using (\ref{J}) and the result of (\ref{GH}).
\be
F_w & = & 1 + F_w (z F_{w-1} + z^2 F_{w-2} \cr
& & + z^2 G_{w-1} (z^{1/2} F_{w-2} + z G_{w-2}) + 2z^{3/2} G_{w-1}). \label{F2}
\ee

Substituting (\ref{G2}) in (\ref{F2}),
\be
F_w F_{w-1} & = & F_{w-1} + z F_w F_{w-1}^2 + z^2 F_w F_{w-1} F_{w-2} + z^2 F_w G_{w-1}^2 \cr
& & + 2 z^{3/2} F_w F_{w-1} G_{w-1} \cr
& = & F_{w-1} +z^2 F_w F_{w-1} F_{w-2} + z F_w F_{w-1} (F_{w-1} + z^{1/2} G_{w-1}) \cr
& & + z^{3/2} F_w G_{w-1} (F_{w-1} + z^{1/2} G_{w-1}) \cr
& = & F_{w-1} +z^2 F_w F_{w-1} F_{w-2} + z^{1/2} G_w F_{w-1} + z G_w G_{w-1}. \label{F3}
\ee

Both (\ref{F2}) and (\ref{F3}) have a term of the form $G_w (z^{1/2} F_{w-1} + z G_{w-1})$. From (\ref{F2}),
\be
z^2 F_w G_{w-1} (z^{1/2} F_{w-2} + z G_{w-2}) & = & F_w - 1 -z F_w F_{w-1}\cr
&- & z^2 F_w F_{w-2} - 2 z^{3/2} F_w G_{w-1}
\ee
\noindent
and from (\ref{F3}),
\be
z^2 F_w G_{w-1} (z^{1/2} F_{w-2} + z G_{w-2}) & = & z^2 F_w (F_{w-1} F_{w-2} \cr
&- &  F_{w-2} - z^2 F_{w-1} F_{w-2} F_{w-3}).
\ee

Equating the two,
\be
F_w & = & 1 + z F_w F_{w-1} + 2 z^{3/2}F_w G_{w-1}  \cr
&+ &   z^2 F_w F_{w-1} F_{w-2} - z^4 F_w F_{w-1} F_{w-2} F_{w-3}. \label{F4}
\ee

Substituting the term $z F_w G_{w-1}$ using (\ref{G2}), we get an expression for $G_w$ in terms of $F_w$'s only.
\be
2 z^{1/2} G_w = F_w - 1 + z F_w F_{w-1} - z^2 F_w F_{w-1} F_{w-2} + z^4 F_w F_{w-1} F_{w-2} F_{w-3} \label{G3}
\ee

Finally, substituting (\ref{G3}) in (\ref{F4}) gives the desired result (\ref{Fw}) $\blacksquare$

\begin{theo}
Let $X_w$ be the generating function for the walk with steps \\ $(1,1),(1,-1),(p,2),(p,-2)$ with $p > 0$. Then $X_w$ satisfies a similar recurrence relation
\be
& X_w = 1- z^{p/2}X_w + (z+z^{p/2})X_w X_{w-1} + (z^{1+p/2}+z^p) X_w X_{w-1} X_{w-2} \cr
& - (z^{3p/2} + z^{2p})X_w X_{w-1} X_{w-2} X_{w-3} + z^{5p/2}X_w X_{w-1} X_{w-2} X_{w-3} X_{w-4}
\ee
\end{theo}

The proof follows exactly the same set of ideas. To start off, we define the same set of generating functions. Equations (\ref{f00}-\ref{f11}) remain the same and equations (\ref{g00}-\ref{g10}) are slightly modified. Following the steps of the previous proof yields the result $\blacksquare$

\section{Numerators and Denominators of $F_w$}
Using (\ref{Fw}), we will now derive a linear recurrence relation for the numerators and denominators of $F_w$. 

\begin{theo}
Let $P_w$ and $\az_w$ be defined as follows. 
\be
& P_0 = 1, & \az_0 =  1 \cr
& P_1 = 1, & \az_1 = 1-z \cr
& P_2 = 1-z, & \az_2 = 1-2z-3z^2 \cr
& P_3 = 1-2z-3z^2, & \az_3 = 1-3z-5z^2-2z^3+z^4 \cr
& P_4 = 1-3z-5z^2-2z^3+z^4, & \az_4 = 1-4z-6z^2+2z^3 \nonumber
\ee

For $w \geq 5$, they are defined recursively by
\be
P_w & = & \az_{w-1} \label{Pw}\\
\az_w & = & (1+z) \az_{w-1}-2z \az_{w-2}- 2z^2 \az_{w-3}\cr
&+& (z^3+z^4) \az_{w-4}-z^5 \az_{w-5} \label{Qw}
\ee

Then, $F_w := \frac{P_w}{\az_w}$ is precisely the generating function for the walk defined earlier satisfying the recurrence relation (\ref{Fw}). 
\end{theo}

These denominators are to basketball what Tchebyshev polynomials are to soccer.

For $w \leq 4$, the generating functions are given by
\be
F_0 & = & 1 \\
F_1 & = & \frac{1}{1-z} \\
F_2 & = & \frac{1-z}{1-2z-3z^2} \\
F_3 & = & \frac{1-2z-3z^2}{1-3z-5z^2-2z^3+z^4} \\
F_4 & = & \frac{1-3z-5z^2-2z^3+z^4}{1-4z-6z^2+2z^3}
\ee
\noindent
and therefore, the initial conditions give the right generating function. To see that (\ref{Pw},\ref{Qw}) imply (\ref{Fw}), divide (\ref{Qw}) by $\az_w$. Then,
\be
1 & = & (1+z) \frac{\az_{w-1}}{\az_w}-2z\frac{\az_{w-2}}{\az_w}- 2z^2\frac{\az_{w-3}}{\az_w}\cr
&+& (z^3+z^4)\frac{\az_{w-4}}{\az_w}-z^5 \frac{\az_{w-5}}{\az_w}.
\ee

But now, using (\ref{Pw})
\be
\frac{\az_{w-1}}{\az_w} & = & F_w, \\
\frac{\az_{w-2}}{\az_w} & = & F_{w-1}F_w, \\
\frac{\az_{w-3}}{\az_w} & = & F_{w-2}F_{w-1}F_w, \\
\frac{\az_{w-4}}{\az_w} & = & F_{w-3}F_{w-2}F_{w-1}F_w, \\
\frac{\az_{w-5}}{\az_w} & = & F_{w-4}F_{w-3}F_{w-2}F_{w-1}F_w.
\ee
which implies (\ref{Fw}) $\blacksquare$

\section{Remarks}
For the sake of completeness, we give references to the number of such basketball games for various values of $w$. For $w=2,\cdots,6$ and $w=\infty$, the sequence of games ending at $n:n$ is in [Sl]. Except for the case of $w=2$, which also arises in some other contexts, all other sequences are new.

Let us now point out why this recurrence is so special! First of all, notice that all terms in (\ref{Fw}) involve only successive generating functions. It is precisely this property that leads to a linear recurrence relation for the denominators. Let us look at this in a little more detail.

Consider the generating functions $F_w, G_w, H_w, J_w$ defined earlier by equations (\ref{F}-\ref{J}). It will not be shown, but it does turn out that the denominators for all four of them are preceisely $\az_w$. Denote their numerators by $P_w,g_w,h_w,j_w$ respectively. Rewriting (\ref{F}-\ref{J}) gives
\be
\frac{P_w}{\az_w} & = & 1 + \frac{P_w}{\az_w \az_{w-1}}(z P_{w-1} + z^{3/2} g_{w-1}+ z^{3/2} h_{w-1} + z^2 j_{w-1}) \\
\frac{g_w}{\az_w} & = & \frac{P_w}{\az_w \az_{w-1}}(z^{1/2} P_{w-1} + z h_{w-1}) \\
\frac{h_w}{\az_w} & = & \frac{P_w}{\az_w \az_{w-1}}(z^{1/2} P_{w-1} + z g_{w-1}) \\
\frac{j_w}{\az_w} & = & \frac{P_{w-1}}{\az_{w-1}} + \frac{g_w}{\az_w \az_{w-1}}(z^{1/2} P_{w-1} + z g_{w-1})
\ee

But notice that $P_w = \az_{w-1}$ and therefore, the first three of these equations are linear but the fourth is not! In fact, if the fourth were also linear, there is no way the recurrence for $P_w$ would terminate uniformly in $w$. The nonlinearity of the fourth equation almost miraculously cancels out excess terms that arise in the fifth order recurrence.

\vspace{1cm}

{\bf REFERENCES}

[BORW] R. Brak, A.L. Owczarek, A. Rechnitzer, S.G. Whittington, {\it A directed walk model of a long chain polymer in a slit with attractive walls}, J. Phys. A, {\bf 38}, 2005, 4309-4325.

[DR] E.A. DiMarzio and R.J. Rubin, {\it Adsorption of a Chain Polymer between Two Plates}, J. Chem. Phys., {\bf 55}, 1971, 4318-36.

[MTW] Keith M. Middlemiss, Glenn M. Torrie and Stuart G. Whittington, {\it Excluded volume effects in the stabilization of colloids by polymers}, J. Chem. Phys., {\bf 66}, 1977, 3227-32.

[R] E.J. Janse van Rensburg, {\it  The statistical mechanics of interacting walks, polygons, animals and vesicles}, Oxford Lecture Series in Mathematics and its Applications, 18. Oxford University Press, Oxford, 2000.

[Sl] N.J.A. Sloane, Sequences A046717,A127617-620,A122951 in the OEIS, \\
\texttt{http://www.research.att.com/$\sim$njas/sequences/Seis.html}

[St] Richard Stanley, Chapter 6 of {\it Enumerative Combinatorics V.2},  Cambridge Studies in Advanced Mathematics, 62. Cambridge University Press, Cambridge, 1999.

[WSCM] Frederick T. Wall, William A. Seitz, John C. Chin and Frederic Mandel, {\it Self-avoiding walks subject to boundary constraints}, J. Chem. Phys., {\bf 67}, 1977, 434-38.

\end{document}